\documentclass[11pt,leqno,twoside]{article}

\usepackage{amsfonts,amsmath,amsthm,amssymb}
\usepackage{enumerate}
\usepackage{graphics,graphicx,subfigure}

\usepackage{color}
\usepackage{todonotes}
\usepackage{cancel}
\usepackage{url}
\usepackage{hyperref}
\usepackage{makeidx}
\usepackage{showidx}
\usepackage{multicol}        % used for the two-column index
\usepackage{xspace}
\usepackage{stmaryrd}        % for \llbracket, \rrbracket
\usepackage{pifont}          % for \ding{227}
\usepackage{fancybox}        % for oval boxes
\usepackage{bm}

 \usepackage{fancyhdr}
\theoremstyle{plain}
 \theoremstyle{definition}
 \newtheorem{lem}{Lemma}
 \newtheorem{defn}[lem]{Definition}
 \newtheorem{thm}[lem]{Theorem}
 \newtheorem{prop}[lem]{Proposition}
 \newtheorem{cor}[lem]{Corollary}
 \newtheorem{notn}[lem]{Notations}
 \newtheorem{pb}[lem]{Problem}
 \newtheorem{form}[lem]{Formulae}
 
 \newtheorem*{rk}{Remark}
 \newtheorem*{com}{Comment}
 \newtheorem*{ex}{Example}
 \theoremstyle{remark}

 \newcommand{\blem}{\begin{lem}}
 \newcommand{\elem}{\end{lem}}
 \newcommand{\bdefn}{\begin{defn}}
 \newcommand{\edefn}{\end{defn}}
 \newcommand{\bthm}{\begin{thm} }
 \newcommand{\ethm}{\end{thm}}
 \newcommand{\bprop}{\begin{prop}}
 \newcommand{\eprop}{\end{prop}}
 \newcommand{\bcor}{\begin{cor}}
 \newcommand{\ecor}{\end{cor}}
 \newcommand{\bnotn}{\begin{notn}}
 \newcommand{\enotn}{\end{notn}}
 \newcommand{\bpb}{\begin{pb}}
 \newcommand{\epb}{\end{pb}}
 \newcommand{\bform}{\begin{form}}
 \newcommand{\eform}{\end{form}}
 \newcommand{\brk}{\begin{rk}}
 \newcommand{\erk}{\end{rk}}
 \newcommand{\bcom}{\begin{com}}
 \newcommand{\ecom}{\end{com}}
 \newcommand{\bex}{\begin{ex}}
 \newcommand{\eex}{\end{ex}}
 \newcommand{\bpf}{\begin{proof}}
 \newcommand{\epf}{\end{proof}}

\newcommand{\cC}{\mathcal{C}}

\newcommand{\cE}{\mathcal{E}}

\newcommand{\cK}{\mathcal{K}}

\newcommand{\cV}{\mathcal{V}}

\newcommand{\bE}{\mathbb{E}}

\newcommand{\bP}{\mathbb{P}}

\newcommand{\bR}{\mathbb{R}}

\newcommand{\be}{\begin{equation}}
\newcommand{\ee}{\end{equation}}
\newcommand{\bal}{\begin{align}}
\newcommand{\eal}{\end{align}}
\newcommand{\ba}{\begin{align*}}
\newcommand{\ea}{\end{align*}}
\newcommand{\bmx}{\begin{matrix}}
\newcommand{\emx}{\end{matrix}}
\newcommand{\bbmx}{\begin{bmatrix}}
\newcommand{\ebmx}{\end{bmatrix}}
\newcommand{\bpmx}{\begin{pmatrix}}
\newcommand{\epmx}{\end{pmatrix}}
\newcommand{\bvmx}{\begin{vmatrix}}
\newcommand{\evmx}{\end{vmatrix}}

\newcommand{\wh}{\widehat}
\newcommand{\wt}{\widetilde}
\newcommand{\f}{\frac}
\newcommand{\df}{\dfrac}

\newcommand{\inc}{\subseteq}

\newcommand{\setm}{\setminus}

\newcommand{\Id}{\mathrm{Id}}

\newcommand{\la}{\lambda}
\newcommand{\La}{\Lambda}
\newcommand{\eps}{\varepsilon}
  
\newcommand{\rev}[1]{{#1}}

\title{\vspace{-20mm}Near-Optimal Estimation of Linear Functionals\\ 
with Log-Concave Observation Errors \medskip\hrule height 1.2pt \vspace{-6mm}}
\author{Simon Foucart\footnote{S. F. is partially supported by grants from the NSF (DMS-2053172) and from the ONR (N00014-20-1-2787).} \, and Grigoris Paouris\footnote{G. P. is partially supported by grants from the NSF (CCF-1900881) and from the Simons Foundation (964286: ``Convexity In High Dimensional Probability'').
Part of this work was carried out while G.P. was a visiting fellow at Princeton University, whose hospitality is greatly appreciated.} \, --- Texas A\&M University}
\date{\vspace{-6mm}\rule{100mm}{0.8pt}}

\newcommand\shorttitle{Near-optimal estimation of linear functionals
with log-concave observation errors}
\newcommand\authors{S. Foucart, G. Paouris}

\begin{document}
\maketitle

%% Add abstract, keywords, and AMS classification
\vspace{-10mm}
\begin{abstract}
This note addresses the question of optimally estimating a linear functional of an object acquired through linear observations corrupted by random noise,
where optimality pertains to a worst-case setting tied to a symmetric, convex, and closed model set containing the object.
It complements the article ``Statistical Estimation and Optimal Recovery'' published in the Annals of Statistics in 1994.
There, Donoho showed (among other things) that, for Gaussian noise,
linear maps provide near-optimal estimation schemes relatively to a performance measure relevant in Statistical Estimation.
Here, we advocate for a different performance measure arguably more relevant in Optimal Recovery.
We show that, relatively to this new measure,
linear maps still provide near-optimal estimation schemes even if the noise is merely log-concave.
Our arguments,
which make a connection to the deterministic noise situation
and bypass properties specific to the Gaussian case,
offer an alternative to parts of Donoho's proof.

\end{abstract}



\vspace{-5mm}
\begin{center}
\rule{100mm}{0.8pt}
\end{center}

%%%%%%%%%%%%%%%%%
%% The main text starts here %%
%%%%%%%%%%%%%%%%%

\section{Introduction}

In this note, we take a second look at the Optimal Recovery problem when random observation errors are present.
As a very brief reminder,
we recall that the Optimal Recovery problem consists in recovering an object $f$---typically a function---from observational data $y_i = \la_i(f)$---typically point evaluations---in a way that is worst-case optimal or near-optimal relatively to a model set~$\cK$.
Here, the difference with this standard scenario is that the observations $y_i$ are corrupted with random additive errors $e_i$, 
so that $y_i = \la_i(f) + e_i$.
Thus, the situation is as follows:
an element $f$ from a Banach space~$F$  is partially known through:
%\vspace{-5mm}
\begin{itemize}
\item some {\sl a priori} information: $f$ belongs to a subset $\cK$ of $F$, i.e.,
$$
f \in \cK,
$$
where $\cK$ is called the model set;
\item some {\sl a posteriori} information: $f$ is inaccurately observed through the actions of some linear functionals $\la_1,\ldots,\la_m \in F^*$, i.e.,
$$
y_i = \la_i(f) + e_i,
\qquad i = 1,\ldots,m.
$$
This is summarized as $y = \La f + e$, where the linear map $\La: F \to \bR^m$ is called the observation map.
Here, $e \in \bR^m$ is a random vector.
\end{itemize}

When estimating $f$, or merely a quantity of interest $Q(f)$ taking values in some Banach space $Z$,
we simply apply a so-called recovery map $\Delta: \bR^m \to Z$ to the available observation vector $y = \La f + e$.
The performance of this recovery map could be assessed, 
for some index $p \in [1,\infty]$, via the global recovery error
\be
\label{GE}
{\rm ge}^{\rm se}_p(\Delta) = \bigg( \sup_{f \in \cK}   \bE \big[ \| Q(f) - \Delta(\La f + e)\|_Z^p  \big] \bigg)^{1/p}.
\ee
We appended a superscript ``${\rm se}$'' because 
this choice is favored in Statistical Estimation,
see e.g. the article~\cite{Don}, which contains the classical result being complemented by this note.
However, we prefer to assess the performance of a recovery map $\Delta: \bR^m \to Z$ via another global recovery error, namely
\be
\label{GE'}
{\rm ge}^{\rm or}_p(\Delta) = \bigg(   \bE \bigg[ \sup_{f \in \cK}  \| Q(f) - \Delta(\La f + e)\|_Z^p  \bigg] \bigg)^{1/p}.
\ee
We appended a superscript ``${\rm or}$'' because we believe
that this choice is better suited to a worst-case perspective,
hence more relevant in Optimal Recovery.
Indeed, suppose that ${\rm ge}^{\rm se}_p(\Delta)$ and ${\rm ge}^{\rm or}_p(\Delta)$ are small,
say bounded by some $\theta$:
Markov's inequality in conjunction with ${\rm ge}^{\rm se}_p(\Delta)^p \le \theta^p$ would naturally yield the statement
$$
\mbox{for all }f \in \cK,
\quad 
\bP \bigg[  \| Q(f) - \Delta(\La f + e)\|_Z  \le \f{\theta}{\eps} \bigg] \ge 1-\eps^p,
$$
while Markov's inequality in conjunction with ${\rm ge}^{\rm or}_p(\Delta)^p \le \theta^p$ would naturally yield the statement
$$ 
\bP \bigg[  \| Q(f) - \Delta(\La f + e)\|_Z  \le \f{\theta}{\eps}  \mbox{ for all }f \in \cK
\bigg] \ge 1-\eps^p.
$$

Of course,  in the absence of observation errors ($e=0$),
these two notions coincide and reduce to a quantity which is independent of $p$,
namely to the global worst-case error (aka distortion)
$$
{\rm gwce}(\Delta) =
\sup_{f \in \cK} \|Q(f) - \Delta(\La f) \|_Z.
$$
In this case, if the model set is symmetric and convex 
(i.e., if $-\cK = \cK$ and $ (1/2) \cK + (1/2) \cK \inc \cK$)
and if $Q: F \to \bR$ is a linear functional,
a classical result of Smolyak  (see \cite{SmoBak} or \cite[Theorem 9.3]{BookDS})
states that ``linear recovery maps are optimal'',
meaning that there exists a linear map $\Delta_{\rm lin}: \bR^m \to \bR$ such that
$$
{\rm gwce}(\Delta_{\rm lin}) 
= \inf_{\Delta: \bR^m \to \bR}
{\rm gwce}(\Delta).
$$ 

In the presence of Gaussian observation errors,
\rev{the previously mentioned seminal work \cite{Don} of Donoho implies that,
although not yielding genuine optimality anymore,  
``linear recovery maps are still near-optimal''.}
Precisely, for $p=1$ and $p=2$,
if $\cK$ is symmetric, convex, closed, and bounded,
if $Q: F \to \bR$ is a linear functional,
and if $e \in \bR^m$ is a mean-zero Gaussian random vector,
then 
there exists a linear map $\Delta_{\rm lin}: \bR^m \to \bR$ such that
$$
{\rm ge}^{\rm se}_p(\Delta_{\rm lin}) 
\le \kappa \times \inf_{\Delta: \bR^m \to \bR}
{\rm ge}^{\rm se}_p(\Delta),
$$ 
where $\kappa$ is an absolute constant not exceeding \rev{$1.23$}.
As a matter of fact,
the validity of this result for $p=1$ implies its validity for all $p \in [1,\infty)$
and even for ${\rm ge}^{\rm or}_p$ \rev{in place of} ${\rm ge}^{\rm se}_p$,  as explained in Subsection~\ref{SubsecCompGEs}.

In this note,
we relax the Gaussianity assumption to the mere requirement that the random vector $e \in \bR^m$ is mean-zero and log-concave.
%and that its covariance matrix is a multiple of the identity,
%i.e., $\bE \big[ e e^\top \big] = \sigma^2 \Id_m$.
Relevant examples include vectors with independent entries distributed according to the Gaussian, Laplace, or uniform distribution.
Uniform distributions on convex sets with appropriately chosen linear structure provide another important example,
see~\cite{BGVV} for some recent developments.
Under the log-concavity assumption,
we still show that ``linear recovery maps are near-optimal'',
but with ${\rm ge}^{\rm or}_p$ \rev{in place of} ${\rm ge}^{\rm se}_p$.
Precisely, we show (Theorem~\ref{ThmMain}) that,
for any $p \ge 1$,
if $\cK$ is symmetric\footnote{Donoho's work drops the assumption that $\cK$ is symmetric and shows that ``affine recovery maps are near-optimal''.
For simplicity of presentation, we did not pursue such a general result.}, convex, and closed (but not necessarily bounded),
if $Q: F \to \bR$ is a linear functional,
and if $e \in \bR^m$ is a mean-zero log-concave random vector,
then 
there exists a linear map $\Delta_{\rm lin}: \bR^m \to \bR$ such that
$$
\rev{
{\rm ge}^{\rm or}_p(\Delta_{\rm lin}) 
\le \kappa_p \times \inf_{\Delta: \bR^m \to \bR}
{\rm ge}^{\rm or}_p(\Delta),
}
$$ 
where $\kappa_p$ is a constant depending only on $p$ that we did not attempt to optimize.

This result is established in Section \ref{SecMain},
where we also point out that our proof supplies streamlined arguments for Donoho's original result from \cite{Don}.
Prior to that, we isolate in Section \ref{SecBG} several ingredients to be relied upon later.
In the spirit of \cite{Don},
we consider  one-dimensional subproblems as a prerequisite for the full problem in
Section \ref{Sec1Dim}, 
where we remark in passing that not all random distributions allow for the near-optimality result.

\section{Background Information}
\label{SecBG}

\subsection{Properties of log-concave random vectors}

%{\color{red} Grigoris, can you check this subsection? 
%In particular, is \cite{MS} the right reference for all this content?}

A \rev{probability} measure $\mu$ on $\bR^m$ is called log-concave if,
for all compact subsets $\cC_0$ and $\cC_1$ of $\bR^m$ and all $\tau \in [0,1]$,
$$
\mu ( (1-\tau) \cC_0 + \tau \cC_1 ) \ge \mu(\cC_0)^{1-\tau} \mu(\cC_1)^\tau.
$$  
A result of C. Borell \rev{ensures} that a log-concave \rev{probability} measure---provided it is not supported on a subspace---satisfies $\mu(\cC) = \int_\cC \pi(x) dx$, $\cC \inc \bR^m$,
for some integrable function $\pi:\bR^m~\to~\bR_+$ such that $-\ln(\pi): \bR^m \to \bR \cup \{\infty\}$ is a convex function. 
\rev{We refer to \cite{BGVV} for the properties of log-concavity listed below.}

A random vector $e \in \bR^m$ is called log-concave if it is distributed according to a probability measure which is log-concave.
The following fact about log-concave random vectors,
known as Borell's lemma,
will be useful later.

\blem
\label{LemBor}
Let $e \in \bR^m$ be a log-concave random vector and let $|\cdot|$ be a seminorm on $\bR^m$.
Then, for any $1 \le p \le q < \infty$,
$$
\big( \bE[ |e|^q ]  \big)^{1/q}
\le C \f{q}{p} \big( \bE[ |e|^p ]  \big)^{1/p},
$$
where the absolute constant $C$ can be taken as $C = {\mathrm e}$.
\elem

Another useful fact for us is that, if $e \in \bR^m$ is a mean-zero log-concave random vector with covariance matrix $\bE [e e^\top] = \sigma^2 \Id_m$ and if $u \in \bR^m$ is an $\ell_2$-normalized vector,
then $\xi = \langle u, e \rangle \in \bR$ is a mean-zero log-concave random variable with variance $\sigma^2$.
We will also rely on the following property of log-concave random variables.
The result is not new, 
but we could not pinpoint the exact statement in the literature.
So, for the reader's convenience,
we provide a proof inspired by an argument of Milman and Pajor from \cite{MP-ini}.
An extension to log-concave random vectors follows from results of \cite[Section 5]{KlCLT}.

\blem
\label{LemLCU}
Let $\pi: \bR \to \bR_+$ be the probability density function of a mean-zero log-concave random variable with variance $\sigma^2$. 
Then 
$$
\pi(x) \ge \f{\delta}{\sigma},
\qquad \mbox{whenever} \qquad
|x| \le \gamma \, \sigma,
$$
where the constants $\delta$ and $\gamma$ can be taken as \rev{$\delta = 1/(2 \sqrt{3} \, {\mathrm e}^2)$} and $\gamma = 1/(5 {\mathrm e}$). 
\elem

\bpf
The argument makes crucial use of a two-sided estimate for $\pi(0)$, namely 
\rev{
\be
\label{EstPi0}
\f{1}{2\sqrt{3} \, {\mathrm e} \,\sigma} \le \pi(0) \le \f{3}{\sigma},
\ee}which goes back to Hensley \cite{H} (see also \cite[Section 2.5]{MP-ini}, \cite{BGVV}, or \rev{\cite[Lemma 2.6]{Kol}}).
\rev{There, the result is stated for a symmetric random variable.
The nonsymmetric case, which follows by applying a result of Fradelizi (\cite[Theorem 2.2.2]{BGVV}),  appears in \cite{BGVV},  see Theorem 2.2.3 and  Lemma 2.2.4.}

\rev{We shall now} prove that $\pi(\pm \gamma \sigma) \ge {\mathrm e}^{-1} \pi(0)$.
\rev{This will imply} that, 
for any $x \in [-\gamma \sigma, \gamma \sigma]$ written as $x = (1-\tau) \times (-\gamma \sigma) + \tau \times (\gamma \sigma)$ for some $\tau \in [0,1]$,
we have \rev{$\pi(x) \ge \pi(-\gamma \sigma)^{1-\tau} \pi(\gamma \sigma)^\tau \ge {\mathrm e}^{-1} \pi(0)$, yielding the announced inequality 
$\pi(x) \ge \delta / \sigma $ with $\delta = 1/(2 \sqrt{3} \, {\mathrm e}^2)$ by invoking \eqref{EstPi0}.}
So let us assume on the contrary that one of $\pi(-\gamma \sigma)$ or $\pi(-\gamma \sigma)$ is smaller than ${\mathrm e}^{-1} \pi(0)$,
e.g. that $\pi(\gamma \sigma) < {\mathrm e}^{-1} \pi(0)$.
Then, for $x \ge \gamma \sigma$,
$$
{\mathrm e}^{-1} \pi(0) > \pi(\gamma \sigma)
= \pi \bigg( \Big(1-\f{\gamma \sigma}{x} \Big) \times 0 + \f{\gamma \sigma}{x} \times x \bigg)
\ge \pi(0)^{1-\gamma \sigma/x} \pi(x)^{\gamma \sigma/x}.
$$
Rearranging the latter, we deduce that $\pi(x) < \pi(0) {\mathrm e}^{-x/(\gamma \sigma)}$ for $x \ge \gamma \sigma$.
It follows that
\rev{
\begin{align*}
\int_{\gamma \sigma}^\infty x \pi(x) dx
& < \pi(0) \int_{\gamma \sigma}^\infty x {\mathrm e}^{-x/(\gamma \sigma)} dx
= \pi(0) \Big[ -\gamma \sigma(x+\gamma \sigma) {\mathrm e}^{-x/(\gamma \sigma)} \Big]_{\gamma \sigma}^\infty
%= \pi(0) \, (\gamma \sigma)^2 \int_1^\infty u {\mathrm e}^{-x/(\gamma \sigma)} du
%= \pi(0) \, (\gamma \sigma)^2 \Big[ -(u+1) {\mathrm e}^{-u} \Big]_1^\infty \\
%&  \le  \f{3}{\sigma} \, (\gamma \sigma)^2 \,  2 {\mathrm e}^{-1}
= \pi(0)\, 2 (\gamma \sigma)^2 {\mathrm e}^{-1}\\
& \le \f{6}{{\mathrm e}} \gamma^2 \sigma,
\end{align*}where the last step utilized \eqref{EstPi0}.}
Moreover, we also have
$$
\int_0^{\gamma \sigma} x \pi(x) dx
\le \gamma \sigma \int_0^{\gamma \sigma} \pi(x) dx
\le \gamma \sigma \int_{-\infty}^\infty \pi(x) dx
=  \gamma \sigma.
$$
Adding these two inequalities, then using the mean-zero property and Lemma \ref{LemBor},
we obtain
$$
\Big(1 + \f{6 \gamma}{{\mathrm e}} \Big) \gamma \sigma
> \int_0^\infty x \pi(x) dx
= \f{1}{2} \int_{-\infty}^{\infty} |x| \pi(x) dx 
\ge \f{1}{2} \, \f{1}{2 {\mathrm e}} \bigg[\int_{-\infty}^\infty x^2 \pi(x) dx \bigg]^{1/2}
= \f{1}{4 {\mathrm e}} \sigma.
$$
We derive the desired contradiction
as soon as $\gamma$ is small enough so that 
$(1 + 6 \gamma/{\mathrm e}) \gamma< 1/(4 {\mathrm e})$,
which occurs with our choice $\gamma = 1/(5 {\mathrm e})$.
\epf

\subsection{Comparison of the two notions of global recovery error}
\label{SubsecCompGEs}

In this subsection,
we compare the notions of global recovery error introduced in \eqref{GE} and \eqref{GE'}.
The results are stated \rev{right below and proved shortly afterwards}.
We note that they are valid when $Q: F \to Z$ is an arbitrary linear map---in particular,  $Q$ need not be a linear functional at this stage and it could even be $Q = \Id_F$.

\bprop
\label{PropCompAll}
\rev{Let $\La: F \to \bR^m$ be an observation map defined on a Banach space $F$
and} let $\Delta: \bR^m \to Z$ be a recovery map for the estimation of a linear map $Q: F \to Z$ \rev{with values in a Banach space $Z$.
Regardless of the model set $\cK$ and the random vector $e \in \bR^m$,
one has,}
for any $1 \le p \le q < \infty$,
$$
{\rm ge}^{\rm se}_p(\Delta) \le {\rm ge}^{\rm se}_q(\Delta),
\qquad \quad
{\rm ge}^{\rm or}_p(\Delta) \le {\rm ge}^{\rm or}_q(\Delta),
\qquad \quad
{\rm ge}^{\rm se}_p(\Delta) \le {\rm ge}^{\rm or}_p(\Delta).
$$ 
\eprop

\bprop
\label{PropCompLin}
\rev{Let $\La: F \to \bR^m$ be an observation map defined on a Banach space $F$
and} let $\Delta_{\rm lin}: \bR^m \to Z$ be a linear recovery map for the estimation of a linear map $Q: F \to Z$ \rev{with values in a Banach space $Z$}.
Given $q \in [1,\infty)$,
if the model set $\cK$ is symmetric and if $e \in \bR^m$ is a log-concave random vector, 
then all the quantities  ${\rm ge}^{\rm se}_p(\Delta_{\rm lin})$ and  ${\rm ge}^{\rm or}_p(\Delta_{\rm lin})$, $1 \le p \le q$,
are comparable up to multiplicative constants that depend only on $q$.
\eprop

Before proving these two statements, we point out a key consequence mentioned in the introduction:
if near-optimality of linear maps is acquired for ${\rm ge}^{\rm se}_1$---as established for Gaussian observation errors in \cite{Don}---then it is automatically acquired for all \rev{${\rm ge}^{\rm se}_p$ and ${\rm ge}^{\rm or}_p$.
Here is a precise statement.

\bcor
Let $\La: F \to \bR^m$ be an observation map defined on a Banach space $F$,
let $\cK$ be a symmetric model set,  and let $e \in \bR^m$ be a log-concave random vector.
When estimating a linear map $Q: F \to Z$ with values in a Banach space $Z$,
if there exist an index $p \in [1,\infty)$ and  a constant $\kappa \ge 1$ such that
$$
\inf_{\substack{\Delta: \bR^m \to \bR \\ \Delta \, {\rm linear}}}
{\rm ge}^{\rm se}_p(\Delta) 
\le \kappa \times \inf_{\Delta: \bR^m \to \bR}
{\rm ge}^{\rm se}_p(\Delta),
$$
then,  for any $q \in [p,\infty)$,  there exists a constant $\kappa_{q} \ge 1$ such that
$$
\inf_{\substack{\Delta: \bR^m \to \bR \\ \Delta \, {\rm linear}}}
{\rm ge}^{\rm se/or}_q(\Delta) 
\le \kappa_{q} \times \inf_{\Delta: \bR^m \to \bR}
{\rm ge}^{\rm se/or}_q(\Delta).
$$
\ecor

\bpf
We just have to consider a linear map $\Delta_{\rm lin}: \bR^m \to Z$ that satisfies ${\rm ge}^{\rm se}_p(\Delta_{\rm lin}) \le \kappa \, {\rm ge}^{\rm se}_p(\Delta)$ for all $\Delta: \bR^m \to Z$
and then write
$$
{\rm ge}^{\rm se/or}_q(\Delta_{\rm lin})
%\underset{{\rm Prop.\ref{PropCompLin}}}{\le}
\le
{\rm const}_q  \, {\rm ge}^{\rm se}_p(\Delta_{\rm lin})
\le {\rm const}_q \,  \kappa \, {\rm ge}^{\rm se}_p(\Delta) 
%\underset{{\rm Prop.\ref{PropCompAll}}}{\le}
\le
{\rm const}_q  \, \kappa \, {\rm ge}^{\rm se/or}_q(\Delta),
$$ 
where the left- and right-most inequalities were due to Propositions \ref{PropCompLin} and \ref{PropCompAll}, respectively.
\epf 

We finish this subsection by providing the missing proofs of the above propositions.}

\bpf[Proof of Propositon \ref{PropCompAll}]
The first two inequalities follow from the fact that,
if $1 \le p \le q < \infty$,
then $\|\cdot\|_{L_p(\mu)} \le \|\cdot\|_{L_q(\mu)} $ for any probability measure $\mu$.
The third inequality is a direct consequence of the general fact that $\sup \, \bE \le \bE \, \sup$.
\epf

\bpf[Proof of Propositon \ref{PropCompLin}]
Fixing a linear recovery map $\Delta_{\rm lin}: \bR^m \to Z$ throughout the proof,
we first claim that it is enough to establish that 
\begin{align}
\label{2Pv1}
\hspace{30mm} {\rm ge}^{\rm or}_q(\Delta_{\rm lin}) & \le C_q \, {\rm ge}^{\rm se}_q(\Delta_{\rm lin}),
&  \hspace{32mm} \hfill \mbox{(log-concavity not required)}\\
\label{2Pv2}
\hspace{30mm} {\rm ge}^{\rm se}_q(\Delta_{\rm lin}) & \le D_q \, {\rm ge}^{\rm se}_1(\Delta_{\rm lin}),
&  \hspace{32mm} \mbox{(log-concavity is required)}
\end{align}
for some constant $C_q,D_q$ depending only on $q$.
Indeed, for $1\le  p \le q$, we would then deduce that
$$
{\rm ge}^{\rm se/or}_p(\Delta_{\rm lin}) 
%\underset{{\rm Prop.}\ref{PropCompAll}}{\le} {\rm ge}^{\rm or}_p(\Delta_{\rm lin}) 
\underset{{\rm Prop.}\ref{PropCompAll}}{\le}  {\rm ge}^{\rm or}_q(\Delta_{\rm lin})
\underset{\eqref{2Pv1}}{\le} C_q \, {\rm ge}^{\rm se}_q(\Delta_{\rm lin})
\underset{\eqref{2Pv2}}{\le} C_q \,  D_q \, {\rm ge}^{\rm se}_1(\Delta_{\rm lin})
%\underset{{\rm Prop.}\ref{PropCompAll}}{\le} C_q \,  D_q \, {\rm ge}^{\rm se}_p(\Delta_{\rm lin})
\underset{{\rm Prop.}\ref{PropCompAll}}{\le} C_q \,  D_q \, {\rm ge}^{\rm se/or}_p(\Delta_{\rm lin}).
$$
In order to establish \eqref{2Pv1} and \eqref{2Pv2},
we now remark that the linearity of $\Delta_{\rm lin}$ allows us to write
\be
\label{GE_lin}
{\rm ge}^{\rm se}_q(\Delta_{\rm lin})^q 
= \sup_{f \in \cK} \bE \Big[ \big\| (Q-\Delta_{\rm lin} \La)f - \Delta_{\rm lin} e \big\|_Z^q  \Big].
\ee
From here, we shall lower-bound this quantity using the symmetry of the model set $\cK$.
\rev{For a fixed~$f$ belonging to $\cK$, since $-f$ also belongs to $\cK$,  we note that}
\begin{align*}
{\rm ge}^{\rm se}_q(\Delta_{\rm lin})^q  
& \ge \max_{\pm} \bE \Big[ \big\| (Q-\Delta_{\rm lin} \La)(\pm f) - \Delta_{\rm lin} e \big\|_Z^q  \Big]\\
%\ge \bE \Big[ \big\| (Q-\Delta_{\rm lin} \La)f - \Delta_{\rm lin} e \big\|_Z^q  \Big]
%= \bE \Big[ \big\| (Q-\Delta_{\rm lin} \La)f + \Delta_{\rm lin} e \big\|_Z^q  \Big]\\
& \ge \f{1}{2} 
\bE \Big[ \big\| (Q-\Delta_{\rm lin} \La)f - \Delta_{\rm lin} e \big\|_Z^q 
+ \big\| -(Q-\Delta_{\rm lin} \La)f - \Delta_{\rm lin} e \big\|_Z^q \Big].
\end{align*}
Using the fact that $a^q + b^q \ge (a+b)^q/2^{q-1}$ for $a,b \ge 0$, it follows that
\begin{align*}
{\rm ge}^{\rm se}_q(\Delta_{\rm lin})^q  
& \ge \f{1}{2^q} 
\bE \Big[ \big( \big\| (Q-\Delta_{\rm lin} \La)f - \Delta_{\rm lin} e \big\|_Z 
+ \big\| -(Q-\Delta_{\rm lin} \La)f - \Delta_{\rm lin} e \big\|_Z \big)^q \Big]\\
& \ge \f{1}{2^q}  \bE \Big[  \max\big\{ 2 \big\| (Q-\Delta_{\rm lin} \La)f \big\|_Z ,
2 \big\|  \Delta_{\rm lin} e \big\|_Z \big\} ^q \Big] = \bE \Big[ \chi_f ^q \Big], 
\end{align*}
where, for later convenience,  we have introduced the random variable
$$
\chi_f = \max\big\{ \big\| (Q-\Delta_{\rm lin} \La)f \big\|_Z ,
\big\|  \Delta_{\rm lin} e \big\|_Z \big\}.
$$
Taking the supremum over $f \in \cK$ now yields the lower bound
$$
{\rm ge}^{\rm se}_q(\Delta_{\rm lin})^q  \ge \sup_{f \in \cK} \bE \Big[ \chi_f ^q \Big].
$$
Turning our attention to ${\rm ge}^{\rm or}_q(\Delta_{\rm lin})$,
the linearity of $\Delta_{\rm lin}$,
a triangle inequality,
and the fact that $(a+b)^q \le 2^{q-1} (a^q + b^q)$ for $a,b \ge 0$
allow us to write
\begin{align}
\nonumber
{\rm ge}^{\rm or}_q(\Delta_{\rm lin})^q 
& = \bE \bigg[ \sup_{f \in \cK}  \big\| (Q-\Delta_{\rm lin} \La)f - \Delta_{\rm lin} e \big\|_Z^q  \bigg]\\
\nonumber
& \le \bE \bigg[ \sup_{f \in \cK}  2^{q-1} \big( \big\| (Q-\Delta_{\rm lin} \La)f \big\|_Z^q + \big\|\Delta_{\rm lin} e \big\|_Z^q  \big) \bigg]\\
\label{LBGELin}
& = 2^{q-1} \bigg( \sup_{f \in \cK}   \big\| (Q-\Delta_{\rm lin} \La)f \big\|_Z^q
+ \bE \Big[ \big\|\Delta_{\rm lin} e \big\|_Z^q \Big]  \bigg).
\end{align}
For any $f \in \cK$,
we have $\big\| (Q-\Delta_{\rm lin} \La)f \big\|_Z \le \chi_f$,
hence $\big\| (Q-\Delta_{\rm lin} \La)f \big\|_Z^q \le \bE \big[ \chi_f^q \big] \le {\rm ge}^{\rm se}_q(\Delta_{\rm lin})^q $,
as well as $\big\|\Delta_{\rm lin} e \big\|_Z \le  \chi_f$,
hence $\bE \big[ \big\|\Delta_{\rm lin} e \big\|_Z^q \big] \le \bE \big[ \chi_f^q \big] \le {\rm ge}^{\rm se}_q(\Delta_{\rm lin})^q $.
This implies that
$$
{\rm ge}^{\rm or}_q(\Delta_{\rm lin})^q  \le 2^q {\rm ge}^{\rm se}_q(\Delta_{\rm lin})^q ,
$$
which is the required inequality \eqref{2Pv1} with $C_q = 2$ (independent of $q$).

For the inequality \eqref{2Pv2},
we come back to \eqref{GE_lin},
use a triangle inequality 
and $(a+b)^q \le 2^{q-1} (a^q + b^q)$ for $a,b \ge 0$
to arrive at
\begin{align*}
{\rm ge}^{\rm se}_q(\Delta_{\rm lin})^q 
& \le \sup_{f \in \cK} \bE \bigg[ 2^{q-1} \big( \big\| (Q-\Delta_{\rm lin} \La)f \big\|_Z^q 
+ \big\| \Delta_{\rm lin} e \big\|_Z^q \big)  \bigg]\\
& = 2^{q-1} \bigg(  \sup_{f \in \cK} \big\| (Q-\Delta_{\rm lin} \La)f \big\|_Z^q 
+ \bE \Big[ \big\| \Delta_{\rm lin} e \big\|_Z^q \Big] \bigg)\\
& \le 2^{q-1} \bigg(  \sup_{f \in \cK} \big\| (Q-\Delta_{\rm lin} \La)f \big\|_Z^q 
+ (C \, q)^q \bE \Big[ \big\| \Delta_{\rm lin} e \big\|_Z \Big]^q \bigg),
\end{align*}
where the last step relied on Borell's lemma (Lemma \ref{LemBor}) for log-concave random vectors.
As before,
for any $f \in \cK$,
we have $\big\| (Q-\Delta_{\rm lin} \La)f \big\|_Z \le \chi_f$,
hence $\big\| (Q-\Delta_{\rm lin} \La)f \big\|_Z \le \bE \big[ \chi_f \big] \le {\rm ge}^{\rm se}_1(\Delta_{\rm lin})$,
as well as $\big\|\Delta_{\rm lin} e \big\|_Z \le  \chi_f$,
hence $\bE \big[ \big\|\Delta_{\rm lin} e \big\|_Z \big] \le \bE \big[ \chi_f \big] \le {\rm ge}^{\rm se}_1(\Delta_{\rm lin})$.
This implies that 
$$
{\rm ge}^{\rm se}_q(\Delta_{\rm lin})^q 
\le 2^{q-1} \Big( {\rm ge}^{\rm se}_1(\Delta_{\rm lin})^q + (C \, q)^q {\rm ge}^{\rm se}_1(\Delta_{\rm lin})^q \Big)
\le 2^q (C \, q)^q {\rm ge}^{\rm se}_1(\Delta_{\rm lin})^q,
$$
which is the required inequality \eqref{2Pv2} with $D_q = 2 \, C \, q$.
\epf

\rev{
\brk
As revealed in the above proof,  
all notions ${\rm ge}_p^{\rm se/or}(\Delta_{\rm lin})$ of global recovery error 
for linear maps $\Delta_{\rm lin}$ are comparable to the maximum (or sum)
of the noiseless global worst-case error ${\rm gwce}(\Delta_{\rm lin}) = \sup_{f \in \cK} \|Q(f) - \Delta_{\rm lin}( \La f) \|_Z$
and the expected noise $\bE [ \|\Delta_{\rm lin} e \|_Z ]$.
\erk
}

\subsection{Optimal estimation with deterministic observation errors}

Throughout this subsection,
it is assumed that the quantity of interest is a linear functional,
in short that $Q \in F^*$.
As for the model set $\cK$, it is assumed to be symmetric and convex,
so it can be thought of in terms of its Minkowski functional (aka gauge function)
$$
|f|_\cK = \inf \{ t>0: f \in t \, \cK\},
\qquad f \in F,
$$
recalling that $|\cdot|_\cK: F \to \bR_+ \cup \{\infty\}$ is a seminorm in the present situation.
Moreover,
we take notice of the equivalence $f \in \cK  \Leftrightarrow  |f|_\cK \le 1 $ when the set $\cK$ is furthermore closed \rev{in the norm of $F$}.

In the accurate setting (where there is no observation \rev{error}),
we have already pointed out that ``linear recovery maps are optimal''.
This remains true in the presence of observation errors modeled deterministically  via the assumption that $e \in \cE$
for some symmetric and convex subset $\cE$ of $\bR^m$.
The relevant example in this note is $\cE = \{ e \in \bR^m: \|e\|_2 \le \sigma\}$,
for which the Minkowski functional is given by $|e|_\cE = \|e\|_2 / \sigma$, $e \in \bR^m$. 
The precise optimality result reads as follows
(\rev{the statements in \eqref{Det2} are} somewhat present in \cite{EF}, but for a specific model set based on approximability).

\bprop
\label{PropDet}
Let $Q: F \to \bR$ be a linear functional.
If the sets $\cK \inc F$ and $\cE \inc \bR^m$ are symmetric, convex, and closed,
then
\begin{align}
\label{Det0}
\inf_{\Delta: \bR^m \to \bR} \, \sup_{f \in \cK, \, e \in \cE} \big| Q(f) - \Delta(\La f+ e) \big|
& = \min_{\Delta_{\rm lin}: \bR^m \to \bR \, {\rm linear}} \, \sup_{ f \in \cK, \, e \in \cE} \big| Q(f) - \Delta_{\rm lin}(\La f+ e) \big|\\
\label{Det1}
& = \min_{a \in \bR^m} \bigg\{ 
\sup_{f \in \cK} \bigg| \Big( Q-\sum_{i=1}^m a_i  \la_i \Big) f \bigg| 
+ \sup_{e \in \cE} \Big| \langle a,e \rangle \Big|  \bigg\}\\
\label{Det2}
& \rev{= \sup_{h \in \cK, \, \La h \in \cE } \big| Q(h) \big| }
= \sup_{h \in F \setm \{0\}} \f{ | Q(h) | }{\max \{ |h|_\cK,  |\La h|_\cE \}} .
\end{align}
\eprop 

\bpf
\rev{The optimality of linear recovery maps expressed by \eqref{Det0} is well known and follows from} a simple reduction\rev{, recalled here,} to the accurate setting.
\rev{Namely}, for any $\Delta: \bR^m \to \bR$,
we interpret the global worst-case error as
$$
\sup_{f \in \cK, \,e \in \cE} |Q(f) - \Delta(\La f+ e) |
= 
\sup_{ (f,e) \in \wt{\cK} } \big| \wt{Q}\big( (f,e) \big) - \Delta\big( \wt{\La}\big( (f, e) \big) \big) \big|,
$$
where the extended quantity of interest $\wt{Q}: F \times \bR^m \to \bR$ is the linear functional defined by $\wt{Q}\big( (f,e) \big)  = Q(f)$
and the extended observation map $\wt{\La}: F \times \bR^m \to \bR^m$ is the linear map defined by $\wt{\La}\big( (f, e) \big) = \La f + e$.
Since the extended model set $\wt{\cK} = \cK \times \cE$ is symmetric and convex,
the classical result of Smolyak about optimality of linear maps applies and justifies the \rev{equality~\eqref{Det0}}. 
The \rev{equality~\eqref{Det1}} is obtained by writing any linear recovery map from $\bR^m$ to $\bR$ as
$\Delta_{\rm lin} = \langle a, \cdot \rangle$ for some $a \in \bR^m$ and by minimizing over $a$ (with some simple manipulations in the mix).
The \rev{first equality in \eqref{Det2}} is also a consequence of Smolyak's result,
since it contains (see e.g. \cite[Theorem~9.3]{BookDS}) the fact that the minimal global worst-case error equals the so-called null error, which is
$$
\rev{
\sup_{ (h,e) \in \wt{\cK}, \, (h,e) \in \ker \wt{\La} } \big| \wt{Q}\big( (h,e) \big) \big|
= \sup_{h \in \cK, \, e \in \cE, \, \La h+ e = 0 } \big| Q(h) \big|
= \sup_{h \in \cK, \, \La h \in \cE } \big| Q(h) \big|.}
$$  
It remains to notice that \rev{$[ h \in \cK \mbox{ and } \La h \in \cE] \Leftrightarrow \max\{ |h|_\cK, | \La h |_\cE \} \le 1$}
and exploit homogeneity to arrive at the \rev{second equality} in \eqref{Det2}.
\epf

\brk
Making sense of \eqref{Det2} implicitly requires that $\max\{ |h|_\cK, | \La h |_\cE \} > 0$ whenever $h \in F \setm \{ 0\}$.
This is actually a common assumption in Optimal Recovery,
at least when $\cE$ is a ball relative to some norm $\|\cdot\|$ on $\bR^m$.
If this assumption was violated,
any recovery map $\Delta: \bR^m \to \bR$ would have infinite global worst-case errors,
so the problem would not even be contemplated in the first place!
Indeed, suppose that we could find a nonzero $h \in F$ such that $\|\La h\| = 0$ and $|h|_\cK=0$,
meaning that $\La h = 0$ and that $(1/t) h \in \cK$ for all $t >0$.
Then, fixing $f_0 \in \cK$ and defining $f_x = f_0 + x h$ for any $x>0$,
we notice that $\La f_x = \La f_0$
and that $f_x \in \cK$---this is because $(1-\eps) f_0 + \eps (x/\eps) h$ belongs to $\cK$ as a convex combination of elements from $\cK$, and hence its limit when $\eps \to 0^+$,
i.e., $f_0 + x h = f_x$, belongs to $\cK$. 
In this case,
the \rev{global recovery error} ${\rm ge}^{\rm or}_1(\Delta)$, say,
cannot be finite independently of $Q \in F^*$,
since 
\begin{align*}
{\rm ge}^{\rm or}_1(\Delta)
& \ge \bE \bigg[ \sup_{x > 0} \big| Q(f_x) - \Delta(\La f_x + e) \big| \bigg]
= \bE \bigg[ \sup_{x > 0} \big| Q(f_0) + x Q(h) - \Delta(\La f_0+ e) \big| \bigg]\\
& \ge \bE \bigg[ \sup_{x>0} \Big( x |Q(h)| - \big| Q(f_0) - \Delta(\La f_0+ e) \big| \Big) \bigg]
= \sup_{x>0} \; x |Q(h)| - \bE \Big[ \big| Q(f_0) - \Delta(\La f_0+ e) \big| \Big].
\end{align*}
The latter is certainly infinite for those linear functionals $Q \in F^*$ such that \rev{$Q(h) \not= 0$}.
\erk

\section{The One-Dimensional Lower Bound}
\label{Sec1Dim}

This section is devoted to the simplest setting of all, namely:
$F = \bR$,  $f \in  \cK = [-\tau, \tau]$,
$m=1$, $y = c f + \xi \in \bR$
with a constant $c \in \bR \setm \{0\}$ and a mean-zero random variable $\xi \in \bR$,
and $Q(f) = b f$ with $b \in \bR$.
The global recovery errors of a map $\Delta: \bR \to \bR$ then reduce, for $1 \le p < \infty$, to
\begin{align}
\label{GE1D}
{\rm ge}^{\rm se}_p(\Delta)^p & = \sup_{f \in [-\tau,\tau]} \bE \bigg[ \big| b f - \Delta(cf+\xi) \big|^p  \bigg],\\
\label{GE'1D}
{\rm ge}^{\rm or}_p(\Delta)^p & = \bE \bigg[  \sup_{f \in [-\tau,\tau]} \big| b f - \Delta(cf+\xi) \big|^p  \bigg] .
\end{align}
When $\xi$ is log-concave,
we shall show that ``linear recovery maps are near-optimal'':
for ${\rm ge}^{\rm or}$,
this is expected since we intend to establish this fact in a more general setting;
for ${\rm ge}^{\rm se}$, it may seem more surprising.
The result for ${\rm ge}^{\rm se/or}_p$ with $p \ge 1$ follows from 
the result for ${\rm ge}^{\rm se}_p$ with $p =1 $, as explained in Subsection \ref{SubsecCompGEs},
and the latter is a consequence of an upper bound for the infimum of ${\rm ge}^{\rm se}_1(\Delta_{\rm lin})$ when $\Delta_{\rm lin} : \bR \to \bR$ is a linear map (Lemma \ref{LemSimpleLin} below)
and of a lower bound for ${\rm ge}^{\rm se}_1(\Delta)$ when $\Delta: \bR \to \bR$ is an arbitrary map
(Lemma \ref{LemSimpleArb} below).
This lower bound is in fact an essential step towards the main result.

\blem
\label{LemSimpleLin}
In the simplest setting,
if $\xi$ is a mean-zero random variable with variance $\sigma^2$, then
$$
\inf_{\Delta_{\rm lin}: \bR \to \bR \, {\rm linear}} {\rm ge}^{\rm se}_1(\Delta_{\rm lin})
\le \inf_{\Delta_{\rm lin}: \bR \to \bR \, {\rm linear}} {\rm ge}^{\rm se}_2(\Delta_{\rm lin})
= \f{|b| \tau \sigma}{ \sqrt{\sigma^2 + c^2 \tau^2}} 
\asymp \f{|b|}{|c|} \min\{ \sigma,  |c| \tau  \}.
$$
\elem

\bpf
The leftmost inequality follows from ${\rm ge}^{\rm se}_1(\Delta_{\rm lin}) \le {\rm ge}^{\rm se}_2(\Delta_{\rm lin})$, see Proposition \ref{PropCompAll}.
The rightmost comparison follows from the two-sided estimate $\max\{ \sigma, |c| \tau \}  \le \sqrt{\sigma^2 + c^2 \tau^2} \le \sqrt{2} \max\{ \sigma, |c| \tau \}$ and some straightforward manipulations.
For the middle equality,
representing the action of $\Delta_{\rm lin}$ as the multiplication by some $a \in \bR$, we have
\begin{align*}
{\rm ge}^{\rm se}_2(\Delta_{\rm lin})^2 & = \sup_{f \in [-\tau,\tau]} \bE \Big[ \big( (b-ac)f - a\xi) \big)^2  \Big]
= \sup_{f \in [-\tau,\tau]} \bE \Big[ ((b-ac)f)^2 - 2 (b-ac)f a\xi + (a\xi)^2  \Big]\\
& = \sup_{f \in [-\tau,\tau]}  \big[ ((b-ac)f)^2 + a^2 \sigma^2 \big]
= (b-ac)^2 \tau^2 + a^2 \sigma^2
= ( c^2 \tau^2 + \sigma^2) a^2 - 2b c \tau^2 a + b^2 \tau^2\\
& = \bigg( \sqrt{\sigma^2 + c^2 \tau^2} a - \f{b c \tau^2}{\sqrt{\sigma^2 + c^2 \tau^2}} \bigg)^2
- \f{b^2 c^2 \tau^4}{\sigma^2 + c^2 \tau^2} + b^2 \tau^2 
\ge \f{b^2 \tau^2 \sigma^2}{\sigma^2 + c^2 \tau^2},
\end{align*}
with equality possible for the choice $a = b c \tau^2/(\sigma^2 + c^2 \tau^2)$.
This justifies the value of the infimum of ${\rm ge}^{\rm se}_2(\Delta_{\rm lin})$ over all linear maps $\Delta_{\rm lin} : \bR \to \bR$.
\epf

\blem
\label{LemSimpleArb}
In the simplest setting,
if $\xi$ is a mean-zero log-concave random variable with variance~$\sigma^2$, then,
for any $\Delta: \bR \to \bR$,
$$
{\rm ge}^{\rm or}_1(\Delta) \ge {\rm ge}^{\rm se}_1(\Delta)
\ge \alpha \f{|b|}{|c|} \min\{ \sigma,  |c| \tau  \},
$$
where the constant $\alpha$ can be taken as $\alpha = \rev{1/ (100 \sqrt{3} {\mathrm e}^4)}$.
\elem

\bpf
Since ${\rm ge}^{\rm or}_1(\Delta) \ge {\rm ge}^{\rm se}_1(\Delta)$ in general,  
see Proposition \ref{PropCompAll}, it suffices to lower-bound $ {\rm ge}^{\rm se}_1(\Delta)$, which takes the form
$$
{\rm ge}^{\rm se}_1(\Delta) = \sup_{f \in [-\tau,\tau]} \int_{-\infty}^\infty |b f - \Delta(c f + x)| \pi(x) dx,
$$
where $\pi$ is the probability density function of the log-concave distribution.
According to Lemma~\ref{LemLCU}, it satisfies $\pi(x) \ge \delta / \sigma$ whenever $|x| \le \gamma \sigma$,
which implies that 
$$
{\rm ge}^{\rm se}_1(\Delta) \ge \sup_{f \in [-\tau,\tau]} \f{\delta}{\sigma} \int_{-\gamma \sigma}^{\gamma \sigma} |b f - \Delta(c f + x)| dx.
$$
Let us introduce the quantity $\nu = \gamma \min\{ \sigma, |c| \tau \}/ (2 |c|)$,
so that $|c| \nu \le  \gamma \sigma / 2$ and $\nu \le \gamma \tau /2 \le \tau$,
ensuring that $\pm \nu \in [-\tau,\tau]$.
We obtain
\begin{align*}
{\rm ge}^{\rm se}_1(\Delta)
& \ge  \f{\delta}{\sigma} \int_{-\gamma \sigma}^{\gamma \sigma} | \pm b \nu - \Delta( \pm c \nu + x)| dx
= \f{\delta}{\sigma} \int_{-\gamma \sigma \pm c \nu}^{\gamma \sigma \pm c \nu} | \pm b \nu - \Delta( y)| dy\\
& \ge \f{\delta}{\sigma} \int_{-\gamma \sigma + |c| \nu}^{\gamma \sigma - |c| \nu} | \pm b \nu - \Delta( y)| dy
\ge \f{\delta}{\sigma} \int_{-\gamma \sigma/2}^{\gamma \sigma/2} | \pm b \nu - \Delta( y)| dy.
\end{align*}
In turn, we deduce that
\begin{align*}
{\rm ge}^{\rm se}_1(\Delta)
& \ge \f{\delta}{\sigma} \int_{-\gamma \sigma/2}^{\gamma \sigma/2} 
\bigg( \f{1}{2} | b \nu - \Delta( y)| + \f{1}{2} | -b \nu - \Delta( y)|  \bigg) dy
\ge \f{\delta}{\sigma} \int_{-\gamma \sigma/2}^{\gamma \sigma/2} |b|\nu dy
=  \f{\delta}{\sigma} \gamma \sigma |b| \nu\\
& = \f{\delta \gamma^2}{2} \f{|b|}{|c|} \min\{ \sigma, |c| \tau \},
\end{align*}
which is the desired inequality with $\alpha = \delta \gamma^2 /2 = \rev{1/ (100 \sqrt{3} {\mathrm e}^4)}$.
\epf

We finish this section by emphasizing that near-optimality of linear recovery maps does not apply 
to all types of distributions for the random observation errors,
even in the simplest setting.
Namely,  we prove that a Rademacher distribution scaled to have variance $\sigma^2$ leads to near-optimality for ${\rm ge}^{\rm se}_2$, say,  if and only if $\sigma \le |c| \tau$,
so near-optimality of linear recovery maps is invalid for large noise level.
This is due to $\inf_{\Delta_{\rm lin}} {\rm ge}^{\rm se}_2 (\Delta_{\rm lin}) \asymp (|b|/|c|) \min \{ \sigma, c |\tau| \}$ (Lemma \ref{LemSimpleLin})
and to the result below.

\bprop
In the simplest setting,
if $\xi$ is the mean-zero random variable with variance $\sigma^2$ defined by $\bP[\xi = - \sigma] = \bP[\xi = + \sigma] = 1/2$,
then
$$
\inf_{\Delta: \bR \to \bR} {\rm ge}^{\rm se/or}_2(\Delta)
\left\{ \bmx
\ge \df{|b|}{|c|} \df{\sigma}{\sqrt{2}} & \mbox{ if } \sigma \le |c| \tau,\\
& \\
= 0 & \mbox{ if } \sigma > |c| \tau.
\emx \right. 
%\quad \mbox{and likewise} \quad
%\inf_{\Delta: \bR \to \bR } {\rm ge}^{\rm or}_2(\Delta)
%\left\{ \bmx
%\ge \df{|b|}{|c|} \df{\sigma}{\sqrt{2}} & \mbox{ if } \sigma \le |c| \tau,\\
% & \\
%= 0 & \mbox{ if } \sigma > |c| \tau.
%\emx \right. 
$$
\eprop

\bpf
{\em Case 1: $\sigma \le |c| \tau$.}
Since ${\rm ge}^{\rm or}_2(\Delta) \ge {\rm ge}^{\rm se}_2(\Delta)$ in general, see Proposition \ref{PropCompAll}, 
we only need to establish the lower bound on ${\rm ge}^{\rm se}_2(\Delta)$ for an arbitrary $\Delta: \bR \to \bR$.
Here,
$$
{\rm ge}^{\rm se}_2(\Delta)^2 = \sup_{f \in [-\tau, \tau]}
\Bigg( \f{1}{2} \big| b f - \Delta(cf+\sigma) \big|^2 +  \f{1}{2} \big| b f - \Delta(cf-\sigma) \big| ^2  \Bigg).
$$
Fixing some $y$ in the interval $[-|c|\tau + \sigma, |c|\tau - \sigma]$,
which is nonempty in this case,  
we consider $f_- = (y-\sigma)/c \in [-\tau,\tau]$, so that $cf_-+\sigma = y$
and let $f_+ = (y+\sigma)/c \in [-\tau,\tau]$, so that $cf_+-\sigma = y$.
We obtain
\begin{align*}
{\rm ge}^{\rm se}_2(\Delta)^2 & \ge 
%\left\{  \bmx
%\f{1}{2} \big| q f_- - \Delta(y) \big|^2 \\ \f{1}{2} \big| q f_+ - \Delta(y) \big|^2
%\emx \right\}
 \max_\pm \f{1}{2} \big| b f_\pm - \Delta(y) \big|^2
\ge  \f{1}{2} \bigg( \f{1}{2} \big| b f_- - \Delta(y) \big|^2 +  \f{1}{2} \big| b f_+ - \Delta(y) \big|^2 \bigg)\\
& \ge \f{1}{8} \Big( \big| b f_- - \Delta(y) \big| +  \big| b f_+ - \Delta(y) \big| \Big)^2
\ge \f{1}{8} \big| b (f_- - f_+) \big|^2\\ 
& = \f{b^2 \sigma^2}{2 c^2}.  
\end{align*}

{\em Case 2: $\sigma > |c| \tau$.}
In view of ${\rm ge}^{\rm se}_2(\Delta) \le {\rm ge}^{\rm or}_2(\Delta)$ again,  we only need to establish that $ {\rm ge}^{\rm or}_2(\Delta) = 0$ for an appropriately chosen recovery map $\Delta: \bR \to \bR$.
This map is defined by
$$
\Delta(y) = \left\{ \bmx
\df{b}{c}(y-\sigma) & \mbox{ if } y>0,\\
& \\
\df{b}{c}(y+\sigma) & \mbox{ if } y<0.
\emx \right.
$$
Keeping in mind that
$$
{\rm ge}^{\rm or}_2(\Delta)^{\rev{2}} = 
\f{1}{2} \sup_{f \in [-\tau, \tau]} \big| b f - \Delta(cf+\sigma) \big|^2 + \f{1}{2}  \sup_{f \in [-\tau, \tau]}  \big| b f - \Delta(cf-\sigma) \big|^2 ,
$$
we notice that, for any $f \in [-\tau,\tau]$, we have $cf + \sigma >0$ and $cf-\sigma < 0$,
so that $\Delta(cf + \sigma ) = b f$ and $\Delta(cf - \sigma ) = b f$.
This immediately implies that ${\rm ge}^{\rm or}_2(\Delta) = 0$.
\epf

\section{The Main Result}
\label{SecMain}

This section finalizes the full justification of this note's message, 
namely that ``linear recovery maps are near-optimal'' for the estimation of a linear functional with log-concave observation errors relatively to the unconventional global recovery error ${\rm ge}^{\rm or}$.
The result is formally stated below.

\bthm
\label{ThmMain}
Let $Q: F \to \bR$ be a linear functional.
If the model set $\cK \inc F$ is symmetric, convex, and closed and if $e \in \bR^m$ is a mean-zero log-concave random vector with invertible covariance matrix,
then there exists a linear map $\Delta_{\rm lin}: \bR^m \to \bR$ such that,
\rev{for any $p \in (1,\infty)$,}
$$
{\rm ge}^{\rm or}_p(\Delta_{\rm lin}) \le \rev{\kappa_p} \times \inf_{\Delta: \bR^m \to \bR} {\rm ge}^{\rm or}_p(\Delta)
$$ 
for some constant \rev{$\kappa_p$} depending only on \rev{$p$}.
\ethm

In what remains, we may and do assume that the invertible covariance matrix is of the form $\bE \big[ e e^\top \big] = \sigma^2 \Id_m$.
Indeed, as a positive definite matrix,
it can be written as $\bE \big[ e e^\top \big] = M M^\top$
for some invertible matrix $M \in \bR^{m \times m}$.
Then we can convert the global errors of a recovery map $\Delta: \bR^m \to \bR$,
given the observation map $\La : F \to \bR^m$,
into the global errors of the recovery map $\wt{\Delta} = \Delta \circ M: \bR^m \to \bR$,
given the observation map $\wt{\La} = M^{-1} \circ \La : F \to \bR^m$,
by virtue of the identity $\Delta(\La f + e) = \wt{\Delta} (\wt{\La} f + \wt{e})$.
Here, $\wt{e} := M^{-1}e \in \bR^m$ is still a mean-zero log-concave random vector (log-concavity is preserved under linear transformations) and, importantly,
its covariance matrix is $\bE \big[ \wt{e} \wt{e}^\top \big] = \rev{M^{-1} \bE \big[ e e^\top \big] (M^{\top})^{-1} = M^{-1} (M M^\top) (M^{\top})^{-1} } = \Id_m$.
Thus, the near-optimality result for the original problem reduces to the near-optimality result for the converted problem,
whose covariance matrix is (a multiple of) the identity.

For $p=1$,  the \rev{result in the latter case} is a consequence of an upper bound for the infimum of ${\rm ge}^{\rm or}_1(\Delta_{\rm lin})$ when $\Delta_{\rm lin} : \bR^m \to \bR$ is a linear map (Lemma \ref{LemLin} below)
and of a lower bound for ${\rm ge}^{\rm or}_1(\Delta)$ when $\Delta: \bR^m \to \bR$ is an arbitrary map
(Lemma \ref{LemArb} below).
For $p \ge 1$, \rev{the result} follows from Propositions~\ref{PropCompAll} and~\ref{PropCompLin}.

\blem
\label{LemLin}
Let $Q: F \to \bR$ be a linear functional.
If the model set $\cK \inc F$ is symmetric, convex, and closed and if $e \in \bR^m$ is a mean-zero random vector with covariance matrix $\bE \big[ e e^\top \big] = \sigma^2 \Id_m$,
then there exists a linear map $\Delta_{\rm lin}: \bR^m \to \bR$ such that 
$$
{\rm ge}^{\rm or}_1(\Delta_{\rm lin}) \le  \sup_{h \in F \setm \{0\}}  \f{|Q(h)|}{ \max\{   \| \La h\|_2/\sigma, |h|_\cK \} }.
$$
\elem

\bpf
Given a linear map $\Delta_{\rm lin}: \bR^m \to \bR$,  according to \eqref{LBGELin}, we have
$$
{\rm ge}^{\rm or}_1(\Delta_{\rm lin})
\le \sup_{f \in \cK} \big| (Q - \Delta_{\rm lin} \La)f \big| + \bE \big[ |\Delta_{\rm lin} e| \big]. 
$$
In view of $\bE [ |\Delta_{\rm lin} e| ] \le (\bE [ |\Delta_{\rm lin} e|^{2} ] )^{1/2}$
and writing $\Delta_{\rm lin} = \langle a, \cdot \rangle$ for some $a \in \bR^m$,
we arrive at
\begin{align*}
{\rm ge}^{\rm or}_1(\Delta_{\rm lin}) & \le  \sup_{f \in \cK} \bigg| \bigg(Q - \sum_{i=1}^m a_i \la_i \bigg)f \bigg| + \Big( \bE \big[ \langle a,  e \rangle^2 \big] \Big)^{1/2}\\
& =  \sup_{f \in \cK} \bigg| \bigg( Q - \sum_{i=1}^m a_i \la_i \bigg)f \bigg| + \sigma \|a\|_2.
\end{align*}
The minimum over $a \in \bR^m$ of the latter coincides with the quantity \eqref{Det1} appearing in Proposition~\ref{PropDet} with $\cE = \{ e \in \bR^m: \|e\|_2 \le \sigma \}$,
i.e.,  with the minimal global worst-case error over $\cK$ and $\cE$.
As such, it also equals \eqref{Det2}.
All in all, we have found a linear map $\Delta_{\rm lin}: \bR^m \to \bR$ such that
$$
{\rm ge}^{\rm or}_1(\Delta_{\rm lin})  \le 
\sup_{h \in F \setm \{0\}}  \f{|Q(h)|}{ \max\{ \| \La h\|_2/\sigma,  |h|_\cK\} },
$$
as desired.
\epf

\blem
\label{LemArb}
Let $Q: F \to \bR$ be a linear functional.
If the model set $\cK \inc F$ is symmetric, convex, and closed and if $e \in \bR^m$ is a mean-zero log-concave random vector with covariance matrix $\bE \big[ e e^\top \big]=\sigma^2 \Id_m$,
then, for any recovery map $\Delta: \bR^m \to \bR$,
$$
{\rm ge}^{\rm or}_1(\Delta) \ge \kappa_1 \times \sup_{h \in F \setm\{0\}}  \f{|Q(h)|}{ \max\{  \| \La h\|_2/\sigma,  |h|_\cK \} },
$$
where the constant $\kappa_1$ can be taken as $\kappa_1 = \rev{1/ (100 \sqrt{3} {\mathrm e}^4)}$.
\elem

\bpf
We start by recalling the expression
$$
{\rm ge}^{\rm or}_1(\Delta) = \bE \bigg[  \sup_{f \in \cK} \big| Q(f) - \Delta(\La f+e) \big|  \bigg].
$$
We decompose $f \in F$ as a (unnormalized) direction $h \in F \setm \{0\}$ and a magnitude $t \in \bR$,
so that $f = t h$.
We set aside the cases $\La h =0$ and $|h|_\cK =0$ for now.
Noticing the equivalence $f \in \cK 
%\Leftrightarrow |t h|_\cK \le 1  
\Leftrightarrow |t| \le 1/|h|_\cK$,
we can write
\begin{align}
\nonumber
{\rm ge}^{\rm or}_1(\Delta)
 & = \bE \bigg[  \sup_{h \in F  \setm \{0\} } \sup_{|t| \le 1/|h|_\cK}  \big| Q(h)t - \Delta((\La h)t+e) \big|  \bigg]\\
 \label{ToExchange}
 & = \bE \bigg[  \sup_{h \in F \setm \{0\} } \sup_{|t| \le 1/|h|_\cK}  
 \Big| Q(h)t - \Delta\Big( \f{\La h}{\|\La h\|_2}( \|\La h\|_2 t + \xi) + e_\perp \Big) \Big|  
 \bigg],
\end{align}
after having decomposed $e \in \bR^m$ as $e = \xi \, \La h/\|\La h\|_2 + e_\perp$,
where $\xi = \langle \La h/\|\La h\|_2, e \rangle \in \bR$ is a mean-zero log-concave random variable with variance $\sigma^2$
and $e_\perp \in \bR^m$ is a random vector orthogonal to~$\La h$.
From $\bE[\sup_{h \in F \setm \{0\} } (\cdot)] \ge \sup_{h \in F \setm \{0\} } \bE[(\cdot)]$, we obtain
${\rm ge}^{\rm or}_1(\Delta) \ge \sup_{h \in F \setm \{0\} } E_h$, where
$$
E_h = \bE_{e_\perp} \bigg[ \bE_\xi \bigg[
\sup_{|t| \le 1/|h|_\cK}  \big| Q(h)t - \wt{\Delta}_{e_\perp} \big( \|\La h\|_2 t + \xi \big) \big|
\, \Big| \, e_\perp \bigg]  \bigg]
$$
 for some appropriately defined map $\wt{\Delta}_{e_\perp}: \bR \to \bR$.
Fixing $e_\perp$,
the inner expectation can be interpreted as the one-dimensional recovery error ${\rm ge}^{\rm or}_1(\wt{\Delta}_{e_\perp})$ given in \eqref{GE'1D}.
Thus, according to Lemma \ref{LemSimpleArb}, it can be lower-bounded as
$$
{\rm ge}^{\rm or}_1(\wt{\Delta}_{e_\perp}) \ge \alpha \,  \f{|Q(h)|}{\| \La h\|_2} \min \Big\{ \sigma, \f{\| \La h \|_2}{|h|_\cK} \Big\}
= \alpha \, \f{|Q(h)|}{ \max\{ \| \La h\|_2/\sigma, |h|_\cK \} }.
$$
This lower bound being independent of $e_\perp$,
it remains a lower bound for $E_h$ itself.
We can therefore conclude that 
$$
{\rm ge}^{\rm or}_1(\Delta) \ge \alpha \, \sup_{h \in F \setm \{ 0\}}  \f{|Q(h)|}{ \max\{ \| \La h\|_2/\sigma, |h|_\cK \} }.
$$
This is, in the generic case,
the desired inequality with $\kappa_1$ equal to the constant $\alpha$ from Lemma \ref{LemSimpleArb}.
It remains to deal with the set-aside cases.
Consider first the case $\La h = 0$, which enforces $|h|_\cK > 0$.
The identity \eqref{ToExchange} is then replaced by 
$$
{\rm ge}^{\rm or}_1(\Delta) = 
\bE \bigg[  \sup_{h \in F \setm \{0\} } \sup_{|t| \le 1/|h|_\cK}
\big| Q(h)t - \Delta  e  \big| \bigg]
\ge \sup_{h \in F \setm \{0\} } \sup_{|t| \le 1/|h|_\cK} |Q(h)| |t|
= \sup_{h \in F \setm \{0\} } \f{|Q(h)|}{|h|_\cK},
$$
where the above inequality used the fact that,
whenever $h \in F \setm \{0\}$ and $|t| \le 1/|h|_\cK$,
the double-supremum is at least $\max_\pm |Q(\pm h) t - \Delta e| = |Q(h)| |t| + |\Delta e| \ge |Q(h)| |t|$.
Thus, the desired inequality is even valid with $\kappa_1 =1$ in this situation.
Consider next the case $|h|_\cK = 0$,
which implies that $f = t h \in \cK$ for any $t \in \bR$ and 
 also enforces $\La h \not= 0$.
 Then the lower bound ${\rm ge}^{\rm or}_1(\Delta) \ge \sup_{h \in F \setm \{0\} } E_h$ still holds with any $\tau >0 $ replacing $1/|h|_\cK$,
 and in particular with $\tau > 0$ large enough so that $\min\{ \sigma, \|\La h\|_2 \, \tau\} = \sigma$.
 Thus, resorting to Lemma \ref{LemSimpleArb} yields $E_h \ge \alpha \, (|Q(h)| / \|\La h\|_2) \sigma$,
 which reduces to the desired inequality with $\kappa_1 = \alpha$ in this situation.
\epf

\brk
As \rev{already} mentioned,
the previous argument is easily adapted to retrieve the result of \cite{Don}
for Gaussian observation errors,
which essentially boils down to establishing the above lower bound for ${\rm ge}^{\rm se}_1(\Delta)$ instead of ${\rm ge}^{\rm or}_1(\Delta)$.
We would first express ${\rm ge}^{\rm se}_1(\Delta)$ as in \eqref{ToExchange}
but with expectation and suprema interchanged.
Then, the benefit of the Gaussian case lies in the independence of $\xi$ and $e_\perp$,
so that we can write
\begin{align*}
{\rm ge}^{\rm se}_1(\Delta) 
& = \sup_{h \in F} \sup_{|t| \le 1/|h|_\cK}  
 \bE_{\xi} \bigg[ \bE_{e_\perp} \bigg[
 \big| Q(h)t - \wt{\Delta}_{e_\perp} \big( \|\La h\|_2 t + \xi \big) \big| \bigg]  \bigg]\\
 & \ge \sup_{h \in F} \sup_{|t| \le 1/|h|_\cK}  
 \bE_{\xi} \bigg[
 \Big| Q(h)t -  \bE_{e_\perp}  \Big[ \wt{\Delta}_{e_\perp} \big( \|\La h\|_2 t + \xi \big) \Big] \Big| \bigg] 
\end{align*}
and invoke the one-dimensional lower bound on ${\rm ge}^{\rm se}_1(\wh{\Delta})$ from Lemma \ref{LemSimpleArb} for the map $\wh{\Delta} = \bE_{e_\perp} \big[ \wt{\Delta}_{e_\perp}(\cdot) \big]$.
\erk

\brk
In closing, we point out that our arguments do not just translate into an existence result.
Indeed, the proof of Lemma \ref{LemLin} reveals that a near-optimal recovery map is provided by a recovery map which is genuinely optimal, albeit with respect to observation errors modeled deterministically via $\cE = \{ e \in \bR^m: \|e\|_2 \le \sigma \}$.
The latter has the form of a linear functional $\langle a^\sharp, \cdot \rangle$,
where $a^\sharp \in \bR^m$ is a minimizer of the convex program \eqref{Det1}.
This program is solvable in many practical situations, including,
as described in \cite{EF}, 
the approximability model sets defined for some finite-dimensional linear subspace $\cV$ of $F$ and some parameter $\eps >0$ by
$$
\cK = \{ f \in F: {\rm dist}_{F}(f, \cV) \le \eps \}.
$$
We remark that this model set is symmetric, convex, and closed, but not bounded.
\erk


\begin{thebibliography}{99}
%\vspace{-4mm}
\bibitem{SmoBak}
N.  S.  Bakhvalov.
{\em On the optimality of linear methods for operator approximation in convex classes of functions. }
USSR Computational Mathematics and Mathematical Physics 11 (1971): 244--249.

\bibitem{BGVV} 
S. Brazitikos, A. Giannopoulos, P. Valettas, and B.-H. Vritsiou. Geometry of Isotropic Convex Bodies.
Mathematical Surveys and Monographs, vol. 196. American Mathematical Society, 2014.

\bibitem{Don}
D. L.  Donoho.
{\em Statistical estimation and optimal recovery.}
The Annals of Statistics 22.1 (1994): 238--270.

\bibitem{EF}
M. Ettehad and S. Foucart. 
{\em Instances of computational optimal recovery: dealing with observation errors.}
SIAM/ASA Journal on Uncertainty Quantification 9.4 (2021): 1438--1456.

 \bibitem{BookDS}
 S. Foucart.
 Mathematical Pictures at a Data Science Exhibition.
 Cambridge University Press, 2022.
 
\bibitem{H} 
D. Hensley.
{\em Slicing convex bodies—bounds for slice area in terms of the body’s covariance.} 
Proceedings of the American Mathematical Society 79.4 (1980): 619--625.
 
\bibitem{KlCLT} 
B. Klartag.
{\em A central limit theorem for convex sets}. 
Inventiones Mathematicae 168.1 (2007): 91--131.

\bibitem{Kol}
A. Koldobsky.
Fourier Analysis in Convex Geometry.
Mathematical Surveys and Monographs, vol. 116. 
American Mathematical Society, 2005.
 
\bibitem{MP-ini} 
V. D. Milman and A. Pajor.
{\em Isotropic position and inertia ellipsoids and zonoids of the unit ball of a normed $n$-dimensional space}.
%In: Lecture Notes in Mathematics, vol. 1376, Springer-Verlag, Berlin,1989, 64--104.  
Geometric aspects of functional analysis (1989): 64--104.
 
% \bibitem{MS}
%V. D. Milman and G. Schechtman.
%Asymptotic Theory of Finite-Dimensional Normed Spaces.
%Lecture Notes in Mathematics, vol. 1200, Springer-Verlag, Berlin, 1986. 
 
% \bibitem{FouLia}
%S. Foucart and C. Liao. 
%{\em Optimal Recovery from Inaccurate Data in Hilbert Spaces: Regularize, but what of the Parameter?}
%Constructive Approximation, to appear.  See also arXiv preprint arXiv:2111.02601 (2021).
%
%\bibitem{Pla}
%L. Plaskota. 
%Noisy Information and Computational Complexity. 
%Cambridge University Press, 1996.

\end{thebibliography}
\end{document}